\documentstyle[12pt]{article}
\def\Bbb R{{\rm \bf R}}
\def\proclaim#1{\vskip2mm{\bf #1}\em}
\def\endproclaim{\em \vskip2mm}
\def\tag#1{\eqno(#1)}
\def\gathered{\begin{array}{c}}
\def\endgathered{\end{array}}
\def\text{\mbox}

\begin{document}

\title {A remark on finding the coefficient
of the dissipative boundary condition
via the enclosure method in the time domain}
\author{Masaru IKEHATA\footnote{
Laboratory of Mathematics,
Institute of Engineering,
Hiroshima University,
Higashi-Hiroshima 739-8527, JAPAN}}
\maketitle

\begin{abstract}
An inverse problem for the wave equation
outside an obstacle with a {\it dissipative boundary
condition} is considered.
The observed data are given by a single
solution of the wave equation generated by an initial data supported on an open ball.
An explicit analytical formula for the computation of the coefficient
at a point on the surface of the obstacle
which is nearest to the center of the support of the initial data is given.

\noindent
AMS: 35R30

\end{abstract}

\section{Introduction}

We consider an inverse obstacle scattering problem which is described by the classical wave equation outside an obstacle
with a dissipative boundary condition.

First we formulate the problem.
Let $D$ be a nonempty bounded open subset of $\Bbb R^3$ with $C^2$-boundary such that
$\Bbb R^3\setminus\overline D$
is connected. Let $\gamma$ be a function belonging to $L^{\infty}(\partial D)$ and satisfy $\gamma\ge 0$.
Let $0<T<\infty$.  Let $B$ be an open ball satisfying $\overline B\cap \overline D=\emptyset$.
We denote by $\chi_B$ the characteristic function of $B$; $p$ and $\eta$ the center and (very small) radius of $B$, respectively.

Let $u=u_B(x,t)$ denote the weak solution of the following initial boundary value problem for the classical wave equation:
$$\displaystyle
\left\{
\begin{array}{ll}
\displaystyle
\partial_t^2u-\triangle u=0 & \text{in}\,(\Bbb R^3\setminus\overline D)\times\,]0,\,T[,\\
\\
\displaystyle
u(x,0)=0 & \text{in}\,\Bbb R^3\setminus\overline D,\\
\\
\displaystyle
\partial_tu(x,0)=\chi_B(x) & \text{in}\,\Bbb R^3\setminus\overline D,\\
\\
\displaystyle
\frac{\partial u}{\partial\nu}-\gamma(x)\partial_tu=0 & \text{on}\,\partial D\times\,]0,\,T[.
\end{array}
\right.
\tag {1.1}
$$
Here $\nu$ denotes the outward normal to $D$ on $\partial D$. The
solution class is taken from \cite{DL} and see also \cite{IEO2}
for its detailed description.

The presence of $\gamma$ affects on the energy of the solution of
(1.1). A formal computation yields
$$\displaystyle
{\cal E}'(t)=-\int_{\partial D}\gamma(x)\vert\partial_t u\vert^2dS\le 0,
$$
where
$$\displaystyle
{\cal E}(t)=\frac{1}{2}\int_{\Bbb R^3\setminus\overline D}(\vert\partial_t u\vert^2+\vert\nabla u\vert^2)dx,\,t\in [0,\,T].
$$
Thus the solution of (1.1) is a model of the wave that loses the energy on the surface of the obstacle.
The distribution of $\gamma$ represents the state of the surface of the obstacle.

In this paper we consider the following problem.

{\bf\noindent Problem.}  Fix a large $T$ (to be determined later).
Assume that both $D$ and $\gamma$
are unknown.  Extract information about the
location and shape of $D$ together with $\gamma$ from the wave field $u_B(x,t)$ given at
all $x\in B$ and $t\in\,]0,\,T[$.

For this problem we have already some solution obtained in \cite{IEO2}.
Note that therein the initial data $\partial_tu(x,0)$ is given by $\chi_B(x)$ multiplied
by a function $f(x)$.  However, we are considering a very small $B$ and so we simplify the problem setting
as above.

Set
$$\displaystyle
w(x)=w_B(x,\tau)=\int_0^Te^{-\tau t}u_B(x,t)dt,\,x\in\Bbb R^3\setminus\overline D,\,\tau>0.
\tag {1.2}
$$

Let $v=v_B(\,\cdot\,,\tau)\in H^1(\Bbb R^3)$ denote the weak solution of the modified Helmholtz equation
$$\displaystyle
(\triangle-\tau^2)v+\chi_B=0\,\,\text{in}\,\Bbb R^3.
\tag {1.3}
$$

We introduce two conditions on $\gamma$:

(A1)  $\exists C'>0\,\,\gamma(x)\le 1-C'$ a.e. $x\in\partial D$;

(A2)  $\exists C'>0\,\,\gamma(x)\ge 1+C'$ a.e. $x\in\partial D$.

Define
$$
\begin{array}{ll}
\displaystyle
I_B(\tau)=\int_B(w-v)dx, & \tau>0.
\end{array}
$$
The function $\tau\longmapsto I_B(\tau)$ is called the {\it indicator function} in the {\it enclosure method}.

The enclosure method goes back to \cite{I1} in which an inverse boundary value problem for the Laplace equation
using a {\it single set} of the Cauchy data has been considered.  Since then the basic idea of the
enclosure method has been developed and applied to several inverse obstacle problems governed
by partial differential equations.  It is a method that finds a domain {\it enclosing} unknown object
from the asymptotic behaviour of the indicator function like above.
In particular,  we have already established the following result.

\proclaim{\noindent Theorem 1.1(\cite{IEO2}).}
If $T$ satisfies
$$\displaystyle
T>2\text{dist}\,(D,B),
\tag {1.4}
$$
then we have:

if (A1) is satisfied, then there exists a $\tau_0>0$ such that, for all $\tau\ge\tau_0$
$$\displaystyle
I_B(\tau)>0;
\tag {1.5}
$$

if (A2) is satisfied, then there exists a $\tau_0>0$ such that, for all $\tau\ge\tau_0$
$$\displaystyle
I_B(\tau)<0.
\tag {1.6}
$$

In both cases the formula
$$
\lim_{\tau\longrightarrow\infty}\frac{1}{2\tau}\log\left\vert I_B(\tau)\right\vert
=-\text{dist}\,(D,B),
\tag {1.7}
$$
is valid.

\endproclaim

Formula (1.7) gives us an information about the geometry
(location) of unknown obstacle. More precisely, recall the {\it
first reflector}
$$\displaystyle
\Lambda_{\partial D}(p)=\{q\in\partial D\,\vert\,\vert q-p\vert
=d_{\partial D}(p)\},
$$
where $d_{\partial D}(p)=\inf_{y\in\partial D}\vert y-p\vert$.
Since $\text{dist}\,(D,B)=d_{\partial D}(p)-\eta$, (1.7) gives us
the largest sphere whose exterior encloses $D$. On the sphere
there exists a point on $\Lambda_{\partial D}(p)$. In addition,
given $\omega\in S^2$ taking a small open ball $B'$ centered at
$p+s\omega$ with a small $s\in]0,\,d_{\partial D}(p)[$, we can
make a decision whether the point $p+d_{\partial D}(p)\omega$
belongs to $\partial D$ or not by sounding $u_{B'}\equiv
u_B\vert_{B=B'}$ on $B'$ over a large but finite time interval via
formula (1.7) with $I_{B'}\equiv I_B\vert_{B=B'}$.

The conclusions (1.5) and (1.6) are byproducts and describe {\it
qualitative property} of  $\gamma$ which affects the state of the
surface of unknown obstacle, roughly speaking, whether $\gamma>>1$
or $\gamma<<1$ in terms of the signature of the value of the
indicator function for large $\tau$.

The proof of Theorem 1.1 is based on the following lower and upper estimates
of the indicator function as $\tau\longrightarrow\infty$ (see \cite{IEO2,ISUR}):

$\bullet$ if $\gamma(x)\ge 0$ for a.e. $x\in\partial D$, then

$$\displaystyle
I_{B}(\tau)\ge \int_{\partial D}
\left(\frac{\partial v}{\partial\nu}-\tau\gamma v\right)v\,dS+O(\tau^{-1}e^{-\tau T});
\tag {1.8}
$$

$\bullet$ if $\gamma(x)\ge C'$ a.e. $x\in\partial D$ for a positive constant $C'$, then

$$\displaystyle
I_{B}(\tau)\le\int_{\partial D}
\left(\frac{\partial v}{\partial\nu}-\tau\gamma v\right)v\,dS
+\frac{1}{\tau}\int_{\partial D}\frac{1}{\gamma}
\left\vert\frac{\partial v}{\partial\nu}-\tau\gamma v\right\vert^2\,dS+O(\tau^{-1}e^{-\tau T}).
\tag {1.9}
$$

However, these inequalities do not yield us the {\it quantitative
information} about the distribution of $\gamma$. The aim of this
paper is to fill this gap and obtain an {\it explicit formula}
which explains the reason for the validity of (1.5) and (1.6)
quantitatively.

In what follows we denote by $B_r(x)$ the open ball centered at $x$ with radius $r$.
To describe the formula it has better to introduce some notion
in differential geometry.
Let $q\in\Lambda_{\partial D}(p)$.
Let $S_q(\partial D)$ and $S_q(\partial
B_{d_{\partial D}(p)}(p))$ denote the {\it shape operators} (or {\it Weingarten maps}) at $q$
of $\partial D$ and $\partial B_{d_{\partial D}(p)}(p)$ with
respect to $\nu_q$ and $-\nu_q$, respectively
(see \cite{O} for the notion of the shape operator).
Since $q$ attains the minimum of the function: $\partial D\ni y\longmapsto
\vert y-p\vert$, we have always $S_q(\partial B_{d_{\partial D}(p)}(p))-S_{q}(\partial D)\ge 0$
as the quadratic form on the common tangent space at $q$.

We introduce the following condition which is stronger than  (A1):

(A1)'  $\exists C'>0\,\,\exists C''>0\,\,\,C''\le \gamma(x)\le 1-C'$ a.e. $x\in\partial D$.

Now we are ready to state the main result in this paper.

\proclaim{\noindent Theorem 1.2.}
Assume that $\partial D$ is $C^3$ and $\gamma\in C^2(\partial D)$.
Assume that $\gamma$ satisfies (A1)' or (A2).
Let $T$ satisfy (1.4).
Assume that the set $\Lambda_{\partial D}(p)$ consists of finite points
and
$$\displaystyle
\text{det}\,(S_q(\partial B_{d_{\partial D}(p)}(p))-S_{q}(\partial D))>0\,\,\,
\forall q\in\Lambda_{\partial D}(p).
\tag {1.10}
$$
Then, we have, as $\tau\longrightarrow\infty$
$$\begin{array}{l}
\displaystyle
\,\,\,\,\,\,
\lim_{\tau\longrightarrow\infty}\tau^4e^{2\tau\text{dist}\,(D,B)}I_B(\tau)
\\
\\
\displaystyle
=
\frac{\pi}{2}
\left(\frac{\eta}{d_{\partial D}(p)}\right)^2
\sum_{q\in\Lambda_{\partial D}(p)}
\frac{1}
{\displaystyle
\sqrt{\text{det}\,(S_q(\partial B_{d_{\partial D}(p)}(p))-S_{q}(\partial D))}
}
\,\frac{1-\gamma(q)}{1+\gamma(q)}.
\end{array}
\tag {1.11}
$$

\endproclaim

Note that, in \cite{IEE} we have considered the case when the boundary condition in (1.1)
is given by the Robin boundary condition
$$\begin{array}{ll}
\displaystyle
\frac{\partial u}{\partial\nu}-\beta(x)u=0 & \text{on}\,\partial D\times\,]0,\,T[,
\end{array}
$$
where $\beta$ is a real-valued function on $\partial D$.
We gave an explicit asymptotic formula of the indicator function in this case.
However, the formula contains: complicated information about the shape of $\partial D$,
more precisely, {\it third-and fourth-order derivatives} of a local representation of $\partial D$
at all the points on $\Lambda_{\partial D}(p)$; the term that contains $\beta$ does not appear
as the leading profile of the indicator function.
These together with (1.11) suggest us that information about the values of $\gamma$ is {\it visible}
rather than those of $\beta$.

We mention here how to make use of Theorem 1.2 in {\it remote sensing}.
Let $p\in\Bbb R^3\setminus\overline D$.  Assume that we have a known point
$q\in\Lambda_{\partial D}(p)$.  How can one find the value of $\gamma$ at $q$ by using
the wave phenomena governed by the wave equation.  From point $p$ let us go a little bit forward to $q$.
We denote by $p'$ the point.  Choose a small open ball $B'$ centered at $p'$ and generate the wave $u_{B'}$.
Compute the indicator function $I_{B'}(\tau)$ by using $v_{B'}$
and $w_{B'}$ via (1.2) by observing $u_{B'}$ on $B'$ over time interval $]0,\,T'[$ for a $T'$
satisfying (1.4) with $B$ replaced with $B'$.
Since the set $\Lambda_{\partial D}(p')$ consists of only the single point $q$
and (1.10) for $p$ replaced with $p'$ is satisfied, applying (1.11)  to the present situation,
one gets the quantity
$$\displaystyle
{\cal F}_{B'}(q)\equiv\frac{\pi}{2}
\left(\frac{\eta'}{d_{\partial D}(p')}\right)^2
\frac{1}
{\displaystyle
\sqrt{\text{det}\,(S_q(\partial B_{d_{\partial D}(p')}(p'))-S_{q}(\partial D))}
}
\,\frac{1-\gamma(q)}{1+\gamma(q)},
$$
where $\eta'$ denotes the radius of $B'$ and we have $d_{\partial D}(p')=d_{\partial D}(p)-\vert p-p'\vert$.
Note also that we have
$$\text{det}\,(S_q(\partial B_{d_{\partial D}(p')}(p'))-S_{q}(\partial D))
=\lambda^2-2H_{\partial D}(q)\lambda+K_{\partial D}(q),
$$
where $\lambda=1/d_{\partial D}(p')$, $H_{\partial D}(q)$ and $K_{\partial D}(q)$ are
the mean and Gauss curvatures at $q$ with respect to $\nu_q$.

Thus if the curvatures $H_{\partial D}(q)$ and $K_{\partial D}(q)$ are known, then from ${\cal F}_{B'}(q)$ one can find
$(1-\gamma(q))/(1+\gamma(q))$ and thus $\gamma(q)$ itself provided $\gamma(q)\not=1$.

By the way, if we do not know the two curvatures at $q$, how can
one find $\gamma(q)$ in the observed wave. One simply way is to
increase the number of observed waves.

Let $p_1$, $p_2$ and $p_3$ are three points on the segment $p\rightarrow q$
different from its end points.
Choose three small balls $B_1$, $B_2$ and $B_3$ with a common radius $\eta'$ centered at $p_1$, $p_2$ and $p_3$, respectively.
Applying the same procedure to $I_{B'}(\tau)$ with $B'=B_j$, $j=1,2,3$ mentioned above,
one gets three quantities:
$$\begin{array}{ll}
\displaystyle
{\cal F}_j\equiv
\frac{2}{\pi}
\left(\frac{d_{\partial D}(p_j)}{\eta'}\right)^2
{\cal F}_{B'}(q)\vert_{B'=B_j}=\frac{A}{\sqrt{\lambda_j^2-2H\lambda_j+K}}, &
\displaystyle
j=1,2,3,
\end{array}
\tag {1.12}
$$
where $A=(1-\gamma(q))/(1+\gamma(q))$, $H=H_{\partial D}(q)$ and $K=K_{\partial D}(q)$;
$\lambda_j=1/d_{\partial D}(p_j)$, $j=1,2,3$.
Then, we have
$$
\displaystyle
\left(
\begin{array}{cc}
\displaystyle
-(\lambda_1{\cal F}_1^2-\lambda_2{\cal F}_2^2)
&
\displaystyle
{\cal F}_1^2-{\cal F}_2^2\\
\\
\displaystyle
-(\lambda_2{\cal F}_2^2-\lambda_3{\cal F}_3^2)
&
\displaystyle
{\cal F}_2^2-{\cal F}_3^2
\end{array}
\right)
\left(\begin{array}{c}
\displaystyle
2 H
\\
\displaystyle
K
\end{array}
\right)
=
\left(\begin{array}{c}
\displaystyle
{\cal F}_2^2\lambda_2^2-{\cal F}_1^2\lambda_1^2
\\
\\
\displaystyle
{\cal F}_3^2\lambda_3^2-
{\cal F}_2^2\lambda_2^2
\end{array}
\right).
\tag {1.13}
$$
Solving this system, we will obtain $H$ and $K$.
Then $A^2$ is uniquely determined by one of three
equations (1.12) or
$$\displaystyle
A^2=\frac{1}{3}
\sum_{j=1}^3{\cal F}_j^2(\lambda_j^2-2H\lambda_j+K).
$$
Since we can know $\gamma(q)>1$ or $\gamma(q)<1$ from one of three
equations (1.12), taking the square root of the both side we
obtain $(1-\gamma(q))/(1+\gamma(q))$ and hence $\gamma(q)$ itself.
Note that the discriminant $M$ for (1.13) has the form
$$\displaystyle
M=(\lambda_3-\lambda_2){\cal F}_3^2{\cal F}_2^2
+(\lambda_2-\lambda_1){\cal F}_2^2{\cal F}_1^2
+(\lambda_1-\lambda_3){\cal F}_1^2{\cal F}_3^2.
$$
However, in general, one can not ensure the non vanishing of $M$ exactly since one can not make
the signature of three numbers $\lambda_3-\lambda_2$, $\lambda_2-\lambda_1$, $\lambda_1-\lambda_3$ the same.

It should be pointed out that there is a well-known result due to Majda \cite{Mo} in
the context of the Lax-Phillips scattering theory for the wave equation.  The boundary condition is the same one
as that of (1.1).  It is assumed that the obstacle is {\it strictly convex}.
He considered the high frequency asymptotics for the {\it scattering amplitude} which can be measured at infinity.
Therefore it is the case when $T=\infty$.  He clarified its leading term as the frequency goes to infinity.
The geometrical information about the obstacle contained in the leading term is only
the Gauss curvature.  Our result contains also the mean curvature as mentioned above.
The information about the coefficient $\gamma$ contained in the formula in the {\it back-scattering case} is essentially same as (1.11).
Note also that if $D$ is {\it convex}, then $\Lambda_{\partial D}(p)$
consists of a single point and (1.10) is satisfied.

Finally we present one simple corollary of Theorem 1.2.
To make the dependence on $\gamma$ clear we denote the indicator function $I_B(\tau)$ by $I_B(\tau;\gamma)$.
Assume that we have $\gamma_1$ and $\gamma_0$ belonging to $C^2(\partial D)$ and
satisfying  (A1)' or (A2).
Under the same assumption on $\partial D$ and $\Lambda_{\partial D}(p)$ in Theorem 1.2,
(1.11) for $\gamma=\gamma_1, \gamma_2$ yields
$$\displaystyle
\lim_{\tau\longrightarrow\infty}
\frac{I_B(\tau;\gamma_1)}
{I_B(\tau;\gamma_0)}
=
\frac{
\displaystyle\sum_{q\in\Lambda_{\partial D}(p)}k_q\frac{1-\gamma_1(q)}{1+\gamma_1(q)}
}
{\displaystyle\sum_{q\in\Lambda_{\partial D}(p)}k_q\frac{1-\gamma_0(q)}{1+\gamma_0(q)}},
\tag {1.14}
$$
where
$$\displaystyle
k_q=\frac{1}
{\displaystyle
\sqrt{\text{det}\,(S_q(\partial B_{d_{\partial D}(p)}(p))-S_{q}(\partial D))}
}.
$$
Then, from the right-hand side on (1.14)  together withe assumption (A1)' or (A2) for $\gamma_0$
we can easily obtain the following estimates and formula.

\proclaim{\noindent Corollary 1.1.}
Assume that $\partial D$ is $C^3$.
Let $\gamma_0$ and $\gamma_1$ belong to $C^2(\partial D)$ and satisfy (A1)' or (A2).
Let $T$ satisfy (1.4).
Assume that the set $\Lambda_{\partial D}(p)$ consists of finite points
and satisfies (1.10).
Then, the limit $\lim_{\tau\longrightarrow\infty} I_B(\tau;\gamma_1)/I_B(\tau;\gamma_0)$ exists
and we have
$$\displaystyle
\min_{q\in\Lambda_{\partial D}(p)}
\frac{
\displaystyle\frac{1-\gamma_1(q)}{1+\gamma_1(q)}
}
{
\displaystyle\frac{1-\gamma_0(q)}{1+\gamma_0(q)}
}
\le
\lim_{\tau\longrightarrow\infty}
\frac{I_B(\tau;\gamma_1)}
{I_B(\tau;\gamma_0)}
\le
\max_{q\in\Lambda_{\partial D}(p)}
\frac{
\displaystyle\frac{1-\gamma_1(q)}{1+\gamma_1(q)}
}
{
\displaystyle\frac{1-\gamma_0(q)}{1+\gamma_0(q)}
}.
\tag {1.15}
$$
In particulr, if $\Lambda_{\partial D}(p)$ consists of a single point $q\in\partial D$, we have
$$\displaystyle
\lim_{\tau\longrightarrow\infty}
\frac{I_B(\tau;\gamma_1)}
{I_B(\tau;\gamma_0)}
=\frac{
\displaystyle\frac{1-\gamma_1(q)}{1+\gamma_1(q)}
}
{
\displaystyle\frac{1-\gamma_0(q)}{1+\gamma_0(q)}
}.
\tag {1.16}
$$

\endproclaim

Estimates on  (1.15) give us some global information about the values of $\gamma_1$ relative to $\gamma_0$ 
at all the points on $\Lambda_{\partial D}(p)$ {\it without} knowing the curvatures.
We need just two observed waves generated by the same initial data on one day which is the case $\gamma=\gamma_0$
and another day the case $\gamma=\gamma_1$.  Formula (1.16) gives a {\it deviation} of the value of $\gamma_1$ from $\gamma_0$
at a {\it monitoring point} on the surface of the obstacle.  Note that given an arbitrary $q\in\partial D$
if $p=q+s\nu_q$ and $s$ is a sufficiently small positive number, then $\Lambda_{\partial D}(p)=\{q\}$ and (1.10) is satisfied.

This paper is organized as follows. In Section 2 we give a proof
of Theorem 1.2. It is based on a rough asymptotic formula of the
indicator function as $\tau\longrightarrow\infty$ which has been
derived in \cite{IEO2}.  The point of the proof of Theorem 1.2 is
to clarify the asymptotic profile of the second term in the
formula. It is stated as Theorem 2.1 and proved in Section 3. In
final section we mention a conclusion together with further
problems. In Appendix we give a proof of Lemma 3.1 which is
essential for that of Theorem 2.1. The proof employs a {\it
reflection argument} developed in \cite{LP} for a characterization
of the right-end point of the {\it scattering kernel}. We have
already used the argument or its modification in the framework of
the enclosure method in the time domain, see \cite{IEE, IMax}.

\section{Proof of Theorem 1.2}

Let $R=w-v$.  The function $R$ satisfies
$$\left\{
\begin{array}{ll}
\displaystyle
(\triangle-\tau^2)R=e^{-\tau T}F  & \text{in}\,\Bbb R^3\setminus\overline D,\\
\\
\displaystyle
\frac{\partial R}{\partial\nu}-\tau\gamma R
=-\left(\frac{\partial v}{\partial\nu}-\tau\gamma v\right)
+e^{-\tau T}G
& \text{on}\,\partial D,
\end{array}
\right.
\tag {2.1}
$$
where
$$
\left\{\begin{array}{l}
\displaystyle
F=F(x,\tau)=\partial_t u(x,T)+\tau u(x,T),\\
\\
\displaystyle
G=G(x)=\gamma(x)u(x,T).
\end{array}
\right.
\tag {2.2}
$$
From \cite{IEO2} we have already known that, as $\tau\longrightarrow\infty$
$$\displaystyle
I_B(\tau)
=J(\tau)+E(\tau)+O(\tau^{-1}e^{-\tau T}),
\tag {2.3}
$$
where
$$\displaystyle
J(\tau)=\int_{\partial D}
\left(\frac{\partial v}{\partial\nu}-\tau\gamma v\right)v\,dS
\tag {2.4}
$$
and
$$\displaystyle
E(\tau)
=\int_{\Bbb R^3\setminus\overline D}(\vert\nabla R\vert^2+\tau^2\vert R\vert^2)dx
+\tau\int_{\partial D}\gamma\vert R\vert^2dS.
\tag {2.5}
$$
See also \cite{ISUR} for a brief explanation about the derivation of  formula (2.3).

Thus the essential part of the proof of Theorem 1.2 should be the study of the asymptotic behaviour of
$J(\tau)$ and $E(\tau)$ as
$\tau\longrightarrow\infty$.  The asymptotic behaviour of $J(\tau)$
can be reduced to study a Laplace-type integral \cite{BH}.
For that of $E(\tau)$ we have the following result which enables us to make a reduction
of the study to a Laplace-type integral.

\proclaim{\noindent Theorem 2.1.}
Assume that $\partial D$ is $C^3$ an $\gamma\in C^2(\partial D)$.
Assume that: $\gamma$ has a positive lower bound; the set $\Lambda_{\partial D}(p)$ consists of finite points
and (1.10) is satisfied; there exists a point $q\in\Lambda_{\partial D}(p)$ such that $\gamma(q)\not=1$.

Let $T$ satisfy
$$\displaystyle
T>\text{dist}\,(D,B).
\tag {2.6}
$$
Then, we have, as $\tau\longrightarrow\infty$
$$
\displaystyle
\lim_{\tau\longrightarrow\infty}\frac{E(\tau)}
{\displaystyle
\int_{\partial D}\left(\frac{\partial v}{\partial\nu}-\tau\gamma v\right)\frac{1-\gamma}{1+\gamma}v \,dS}=1.
\tag {2.7}
$$
\endproclaim

The proof of Theorem 2.1 is given in Section 3. We continue to
proceed the proof of Theorem 1.2.

It is well known that the Laplace method under the assumption
that $\Lambda_{\partial D}(p)$ is finite and satisfies (1.10),
we have
$$\begin{array}{l}
\displaystyle
\,\,\,\,\,\,
\lim_{\tau\longrightarrow\infty}
\tau  e^{2\tau d_{\partial D}(p)}
\int_{\partial D}A(x)
\frac{e^{-2\tau\vert x-p\vert}}
{\vert x-p\vert^2}dS
\\
\\
\displaystyle
=\frac{\pi}{d_{\partial D}(p)^2}
\sum_{q\in\Lambda_{\partial D}(p)}
\frac{A(q)}
{\displaystyle
\sqrt{\text{det}\,(S_q(\partial B_{d_{\partial D}(p)}(p))-S_q(\partial D))}},
\end{array}
\tag {2.8}
$$
where $A\in C^1(\partial D)$.  See \cite{BH}, for example.  The point is that the Heassian of the function
$\partial D\ni x\longmapsto \vert x-p\vert$ at $q\in\Lambda_{\partial D}(p)$ is given by
the operator $S_q(\partial B_{d_{\partial D}(p)}(p))-S_q(\partial D)$.  See, for example, \cite{IEE} for this point.

The weak solution $v$ of (1.3) is explicitly given by the formula
$$\displaystyle
v(x)=\frac{1}{4\pi}\int_B\frac{e^{-\tau\vert x-y\vert}}{\vert x-y\vert}dy.
$$
By the mean value theorem for the modified Helmholtz equation
\cite{CH}, we have, for all $x\in\Bbb R^3\setminus\overline B$
$$\displaystyle
\frac{1}{4\pi}
\int_B\frac{e^{-\tau\vert x-y\vert}}{\vert x-y\vert}
dy
=
\frac{\varphi(\tau\eta)}
{\tau^3}
\frac{
e^{-\tau\vert x-p\vert}}
{\vert x-p\vert},
$$
where $\varphi(\xi)=\xi\cosh\,\xi-\sinh\,\xi$.
Therefore $v$ outside $B$ has the explicit form
$$\displaystyle
v(x)=\frac{\varphi(\tau\eta)}
{\tau^3}
\frac{
e^{-\tau\vert x-p\vert}}
{\vert x-p\vert}
\tag {2.9}
$$
and note that
$$\displaystyle
\varphi(\tau\eta)
=\frac{\tau\eta e^{\tau\eta}}{2}(1+O(\tau^{-1})).
\tag {2.10}
$$

Set
$$\displaystyle
\tilde{v}(x)=\frac{e^{-\tau\vert x-p\vert}}{\vert x-p\vert}.
$$
Let $x\in\partial D$.
We have
$$
\begin{array}{l}
\displaystyle
\frac{\partial\tilde{v}}{\partial\nu}
=\left(\tau+\frac{1}{\vert x-p\vert}\right)\frac{p-x}{\vert x-p\vert}\cdot\nu\,\tilde{v}
\end{array}
\tag {2.11}
$$
and thus
$$\displaystyle
\frac{\partial\tilde{v}}{\partial\nu}
-\tau\gamma\tilde{v}
=
\tilde{v}
\left\{
\tau\left(\frac{p-x}{\vert x-p\vert}\cdot\nu-\gamma\right)
+\frac{1}{\vert x-p\vert}\frac{p-x}{\vert x-p\vert}\cdot\nu\right\}.
\tag {2.12}
$$
Then, it follows from (2.8) that
$$\begin{array}{l}
\displaystyle
\,\,\,\,\,\,
\lim_{\tau\longrightarrow\infty}
e^{2\tau\,d_{\partial D}(p)}\int_{\partial D}\left(\frac{\partial\tilde{v}}{\partial\nu}
-\tau\gamma\tilde{v}\right)\tilde{v}dS
\\
\\
\displaystyle
=\frac{\pi}{d_{\partial D}(p)^2}
\sum_{q\in\Lambda_{\partial D}(p)}\frac{1-\gamma(q)}
{\displaystyle
\sqrt{\text{det}\,(S_q(\partial B_{d_{\partial D}(p)}(p))-S_q(\partial D))}}.
\end{array}
\tag {2.13}
$$
Note that we have made use of the fact that, for all $x\in\,\Lambda_{\partial D}(p)$ the unit
vector $(p-x)/\vert x-p\vert$ coincides with $\nu$ at $x$.

Since $v$ has the form (2.9), (2.4) gives
$$\displaystyle
\int_{\partial D}\left(\frac{\partial\tilde{v}}{\partial\nu}
-\tau\gamma\tilde{v}\right)\tilde{v}dS
=\left(\frac{\tau^3}{\varphi(\tau\eta)}\right)^2J(\tau).
$$
Substituting this into (2.13) and using (2.10), we obtain
$$\begin{array}{l}
\displaystyle
\,\,\,\,\,\,
\lim_{\tau\longrightarrow\infty}
\tau^4 e^{2\tau\text{dist}\,(D,B)}J(\tau)
\\
\\
\displaystyle
=\frac{\pi}{4}\left(\frac{\eta}{d_{\partial D}(p)}\right)^2
\sum_{q\in\Lambda_{\partial D}(p)}\frac{1-\gamma(q)}
{\displaystyle
\sqrt{\text{det}\,(S_q(\partial B_{d_{\partial D}(p)}(p))-S_q(\partial D))}}.
\end{array}
\tag {2.14}
$$
Similarly, we obtain
$$
\begin{array}{l}
\displaystyle
\,\,\,\,\,\,
\lim_{\tau\longrightarrow\infty}
\tau^4 e^{2\tau\text{dist}\,(D,B)}\int_{\partial D}\left(\frac{\partial v}{\partial\nu}-\tau\gamma v\right)\frac{1-\gamma}{1+\gamma}v \,dS
\\
\\
\displaystyle
=\frac{\pi}{4}\left(\frac{\eta}{d_{\partial D}(p)}\right)^2
\sum_{q\in\Lambda_{\partial D}(p)}\frac{1}
{\displaystyle
\sqrt{\text{det}\,(S_q(\partial B_{d_{\partial D}(p)}(p))-S_q(\partial D))}}
\frac{(1-\gamma(q))^2}{1+\gamma(q)}
\end{array}
\tag {2.15}
$$
provided $T$ satisfies (2.6).  Note that this yields also that,
for sufficiently large $\tau$ the denominator of (2.7) is positive
under the condition $\gamma(q)\not=1$ for a $q\in\Lambda_{\partial
D}(p)$.

Since we have
$$\displaystyle
\tau^4e^{2\tau\text{dist}\,(D,B)}O(\tau^{-1}e^{-\tau T})=O(\tau^3 e^{-\tau(T-2\text{dist}\,(D,B))})
$$
and
$$\displaystyle
1-\gamma(q)+\frac{(1-\gamma(q))^2}{1+\gamma(q)}
=2\frac{1-\gamma(q)}{1+\gamma(q)},
$$
from (2.3), (2.14), (2.7) and (2.15), we obtain (1.11) provided $T$ satisfies (1.4).

{\bf\noindent Remark 2.1.} It should be emphasized that formula
(2.7) is valid for $T$ satisfying (2.6) which includes the case
when (1.4) is not satisfied.  However, this is reasonable.
Formula (2.7) is concerned with the energy of the {\it reflected}
solution arrived at the surface of the obstacle, thus something
about the {\it incident} field should be encoded therein if $T$ is
greater than the first arriving time $\text{dist}\,(D,B)$ at the
surface of the obstacle.  Formula (2.7) clarifies its {\it
principle term} explicitly.

{\bf\noindent Remark 2.2.}
Theorem 2.1 suggests us inequalities (1.8) and (1.9) are best possible in some sense.
The reason for this is the following.
From (2.12) for $\gamma\equiv 1$, roughly speaking, we have
$$\begin{array}{ll}
\displaystyle
\frac{\partial v}{\partial\nu}\sim \tau v & \text{on}\,\,\Lambda_{\partial D}(p)
\end{array}
$$
and thus we can expect
$$\displaystyle
J(\tau)\sim\tau\int_{\partial D}(1-\gamma)v^2\,dS
$$
and also
$$\displaystyle
\frac{1}{\tau}\int_{\partial D}\frac{1}{\gamma}
\left\vert\frac{\partial v}{\partial\nu}-\tau\gamma v\right\vert^2\,dS
\sim
\tau\int_{\partial D}\frac{1-\gamma}{\gamma}v^2\,dS.
$$
Hence we can expect
$$\displaystyle
J(\tau)+\frac{1}{\tau}\int_{\partial D}\frac{1}{\gamma}
\left\vert\frac{\partial v}{\partial\nu}-\tau\gamma v\right\vert^2\,dS
\sim
\tau\int_{\partial D}\left\{(1-\gamma)+\frac{(1-\gamma)^2}{\gamma}\right\}v^2\,dS.
$$
On the other hand, we can expect also
$$\displaystyle
\int_{\partial D}\left(\frac{\partial v}{\partial\nu}-\tau\gamma v\right)\frac{1-\gamma}{1+\gamma}v \,dS
\sim
\tau\int_{\partial D}\frac{(1-\gamma)^2}{1+\gamma}v^2dS.
$$
Thus from (2.3), (2.4) and (2.7) we can expect
$$\displaystyle
I_B(\tau)
\sim
\tau\int_{\partial D}
\left\{
(1-\gamma)+
\frac{(1-\gamma)^2}{1+\gamma}\right\}v^2\,dS.
$$
Then, we see that (1.8) and (1.9) correspond to the following
trivial inequalities
$$\displaystyle
1-\gamma\le
 (1-\gamma)+\frac{(1-\gamma)^2}{1+\gamma}
\le
(1-\gamma)+\frac{(1-\gamma)^2}{\gamma}.
$$

\section{Proof of Theorem 2.1}

We denote by $x^r$ the reflection across $\partial D$ of the point $x\in\Bbb R^3\setminus D$ with
$d_{\partial D}(x)<2\delta_0$ for a sufficiently small $\delta_0>0$.
It is given by $x^r=2q(x)-x$, where $q(x)$ denotes the unique point on $\partial D$ such that
$d_{\partial D}(x)=\vert x-q(x)\vert$.  Note that $q(x)$ is $C^2$ for $x\in\Bbb R^3\setminus D$ with $d_{\partial D}(x)<2\delta_0$ if $\partial D$ is $C^3$ (see \cite{GT}).
Define $\tilde{\gamma}(x)=\gamma(q(x))$ for $x\in\Bbb R^3\setminus D$ with $d_{\partial D}(x)<2\delta_0$.
The function $\tilde{\gamma}$ is $C^2$ therein and coincides with $\gamma(x)$ for $x\in\partial D$.
It is easy to see that $\partial \tilde{\gamma}/\partial\nu=0$ on $\partial D$.

Choose a cutoff function $\phi_{\delta}\in C^{2}(\Bbb R^3)$ with
$0<\delta<\delta_0$ which satisfies
$0\le\phi_{\delta}(x)\le 1$; $\phi_{\delta}(x)=1$ if $d_{\partial D}(x)<\delta$;
$\phi_{\delta}(x)=0$ if $d_{\partial D}(x)>2\delta$;
$\vert\nabla\phi_{\delta}(x)\vert\le C\delta^{-1}$;
$\vert\nabla^2\phi_{\delta}(x)\vert\le C\delta^{-2}$.

Define
$$\displaystyle
R_0(x)=\frac{1-\tilde{\gamma}(x)}{1+\tilde{\gamma}(x)}\phi_{\delta}(x)v(x^r).
$$
Since we have
$$\begin{array}{ll}
\displaystyle
v(x^r)=v(x),\,\,\frac{\partial}{\partial\nu}v(x^r)=-\frac{\partial v}{\partial\nu}(x)
& \text{on}\,\partial D,
\end{array}
$$
we obtain
$$
\begin{array}{ll}
\displaystyle
\frac{\partial R_0}{\partial\nu}
&
\displaystyle
=
\frac{\partial }{\partial\nu}
\left(
\frac{1-\tilde{\gamma}}{1+\tilde{\gamma}}
\right)v
-\frac{1-\gamma}{1+\gamma}
\frac{\partial v}{\partial\nu}
\\
\\
\displaystyle
&
\displaystyle
=-\frac{1-\gamma}{1+\gamma}
\frac{\partial v}{\partial\nu}
\end{array}
$$
and thus
$$\displaystyle
\frac{\partial R_0}{\partial\nu}
-\tau\gamma R_0
=
-\frac{1-\gamma}{1+\gamma}
\frac{\partial v}{\partial\nu}
-\tau\gamma
\frac{1-\gamma}{1+\gamma}v.
\tag {3.1}
$$
Now define
$$\displaystyle
R_1=R-R_0.
$$
A combination of (3.1) and the boundary condition in (2.1) gives
$$\begin{array}{l}
\displaystyle
\,\,\,\,\,\,
\frac{\partial R_1}{\partial\nu}-\tau\gamma R_1\\
\\
\displaystyle
=
-\left(\frac{\partial v}{\partial\nu}-\tau\gamma v\right)
+\frac{1-\gamma}{1+\gamma}
\frac{\partial v}{\partial\nu}
+\tau\gamma\frac{1-\gamma}{1+\gamma}v
+e^{-\tau T}G\\
\\
\displaystyle
=-\frac{2\gamma}{1+\gamma}
\left(\frac{\partial v}{\partial\nu}-\tau v\right)
+e^{-\tau T}G,
\end{array}
$$
that is
$$
\displaystyle
\frac{\partial R_1}{\partial\nu}-\tau\gamma R_1
=-\frac{2\gamma}{1+\gamma}
\left(\frac{\partial v}{\partial\nu}-\tau v\right)
+e^{-\tau T}G.
\tag {3.2}
$$
Define
$$\displaystyle
v_1=\frac{2\gamma}{1+\gamma}
\left(\frac{\partial v}{\partial\nu}-\tau v\right).
\tag {3.3}
$$
It follows from (2.1)  and (3.2) that $R_1$ satisfies
$$\left\{
\begin{array}{ll}
\displaystyle
(\triangle-\tau^2)R_1=-(\triangle-\tau^2)R_0+e^{-\tau T}F   & \text{in}\,\Bbb R^3\setminus\overline D,\\
\\
\displaystyle
\frac{\partial R_1}{\partial\nu}-\tau\gamma R_1
=-v_1+e^{-\tau T}G  & \text{on}\,\partial D.
\end{array}
\right.
\tag {3.4}
$$
It follows from (2.1) and (2.5) that
$$\displaystyle
E(\tau)
=\int_{\partial D}\left(\frac{\partial v}{\partial\nu}-\tau\gamma v\right)\,R dS
-e^{-\tau T}
\left(\int_{\partial D}GR dS+\int_{\Bbb R^3\setminus\overline D} FRdx\right)
$$
and thus
$$\begin{array}{ll}
\displaystyle
E(\tau)
=
&
\displaystyle
\int_{\partial D}\left(\frac{\partial v}{\partial\nu}-\tau\gamma v\right)\frac{1-\gamma}{1+\gamma}v \,dS\\
\\
&
\displaystyle
+\int_{\partial D}\left(\frac{\partial v}{\partial\nu}-\tau\gamma v\right)\,R_1 dS
-e^{-\tau T}
\left(\int_{\partial D}GR dS+\int_{\Bbb R^3\setminus\overline D} FRdx\right).
\end{array}
\tag {3.5}
$$

In Lemma 2.1 and (2.28) in \cite{IEO2} we have already known that, as $\tau\longrightarrow\infty$
$$\left\{
\begin{array}
{l}
\displaystyle
\Vert R\Vert_{L^2(\Bbb R^3\setminus\overline D)}=
O(e^{-\tau\,\text{dist}\,(D,B)}+e^{-\tau T}),
\\
\\
\displaystyle
\Vert \nabla R\Vert_{L^2(\Bbb R^3\setminus\overline D)}=
O(\tau(e^{-\tau\,\text{dist}\,(D,B)}+e^{-\tau T})),
\\
\\
\displaystyle
\Vert R\Vert_{L^2(\partial D)}=
O(\tau^{1/2}(e^{-\tau\,\text{dist}\,(D,B)}+e^{-\tau T})).
\end{array}
\right.
$$
Noting also the form (2.2) of $F$ and $G$ and the solution class of (1.1) (see \cite{IEO2}), we see that (3.5) becomes
$$\begin{array}{ll}
\displaystyle
E(\tau) =
&
\displaystyle
\int_{\partial D}\left(\frac{\partial v}{\partial\nu}-\tau\gamma v\right)\frac{1-\gamma}{1+\gamma}v \,dS\\
\\
\displaystyle
&
\displaystyle
+\int_{\partial D}\left(\frac{\partial v}{\partial\nu}-\tau\gamma v\right)\,R_1 dS
+\tau e^{-\tau T}(e^{-\tau\,\text{dist}\,(D,B)}+e^{-\tau T}).
\end{array}
\tag {3.6}
$$

Now we study the asymptotic behaviour of $\Vert R_1\Vert_{L^2(\partial D)}$ as $\tau\longrightarrow\infty$.
Multiplying both sides of the first equation on (3.4) with $R_1$ and applying integration by parts,
we obtain
$$\begin{array}{l}
\displaystyle
\,\,\,\,\,\,
\int_{\Bbb R^3\setminus\overline D}
(\vert\nabla R_1\vert^2+\tau^2\vert R_1\vert^2)dx
+\tau\int_{\partial D}\gamma\vert R_1\vert^2 dS\\
\\
\displaystyle
\,\,\,
+e^{-\tau T}
\left(\int_{\partial D} GR_1dS
+\int_{\Bbb R^3\setminus\overline D}FR_1dS\right)
\\
\\
\displaystyle
=\int_{\partial D}v_1R_1dS+
\int_{\Bbb R^3\setminus\overline D}(\triangle-\tau^2)R_0\cdot R_1 dx.
\end{array}
$$
Rewrite this as
$$\begin{array}{l}
\displaystyle
\int_{\Bbb R^3\setminus\overline D}
\left(2\vert\nabla R_1\vert^2+2\tau^2\left\vert R_1+\frac{e^{-\tau T}F}{2\tau^2}\right\vert^2\right)\,dx
\\
\\
\displaystyle
\,\,\,
+2\tau\int_{\partial D}\gamma
\left\vert R_1-\frac{1}{2\tau\gamma}
(v_1-e^{-\tau T}G)
\right\vert^2\,dS\\
\\
\displaystyle
=\frac{e^{-2\tau T}}{2\tau^2}
\int_{\Bbb R^3\setminus\overline D}\vert F\vert^2\,dx
+\frac{1}{2\tau}
\int_{\partial D}\frac{1}{\gamma}\vert v_1-e^{-\tau T}G\vert^2\,dS
\\
\\
\displaystyle
\,\,\,
+2\int_{\Bbb R^3\setminus\overline D}(\triangle-\tau^2)R_0\cdot R_1 dx.
\end{array}
\tag {3.7}
$$

Since we have
$$
\left\{
\begin{array}{l}
\displaystyle
\tau^2\vert R_1\vert^2
\le
2\tau^2\left\vert R_1+\frac{e^{-\tau T}}{2\tau^2}\right\vert^2
+\frac{e^{-2\tau T}\vert F\vert^2}{2\tau^2},\\
\\
\displaystyle
\tau\gamma\vert R_1\vert^2
\le
2\tau\gamma\left\vert R_1-\frac{1}{2\tau\gamma}(v_1-e^{-\tau T}G)\right\vert^2
+\frac{1}{2\tau\gamma}
\vert v_1-e^{-\tau T}G\vert^2,
\end{array}
\right.
$$
it follows from (3.7) that
$$\begin{array}{ll}
\displaystyle
E_1(\tau)
&
\displaystyle
\le
\frac{e^{-2\tau T}}{\tau^2}
\int_{\Bbb R^3\setminus\overline D}\vert F\vert^2 dx
+
\frac{1}{\tau}\int_{\partial D}\frac{1}{\gamma}
\vert v_1-e^{-\tau T}G\vert^2\,dx
\\
\\
\displaystyle
&
\displaystyle
\,\,\,
+2\int_{\Bbb R^3\setminus\overline D}(\triangle-\tau^2)R_0\cdot R_1 dx,
\end{array}
\tag {3.8}
$$
where
$$\displaystyle
E_1(\tau)
=\int_{\Bbb R^3\setminus\overline D}
(\vert\nabla R_1\vert^2+\tau^2\vert R_1\vert^2)\,dx
+\tau\int_{\partial D}\gamma\vert R_1\vert^2\,dS.
$$
Here we make use of a rough estimate:
$$\displaystyle
\vert v_1-e^{-\tau T}G\vert^2
\le 2(\vert v_1\vert^2+e^{-2\tau T}\vert G\vert^2).
$$
Then (3.8) gives
$$\displaystyle
E_1(\tau)
\le
\frac{1}{\tau}\int_{\partial D}\frac{2}{\gamma}
\vert v_1\vert^2\,dx
+2\int_{\Bbb R^3\setminus\overline D}(\triangle-\tau^2)R_0\cdot R_1 dx
+O(e^{-2\tau T}).
\tag {3.9}
$$
Here we claim

\proclaim{\noindent Lemma 3.1.}
If $\delta=\tau^{-1/2}$, then we have
$$\displaystyle
\left\vert\int_{\Bbb R^3\setminus\overline D}(\triangle-\tau^2)R_0\cdot R_1 dx\right\vert
\le C\tau^{-1/2}J_0(\tau)^{1/2}
E_1(\tau)^{1/2},
\tag {3.10}
$$
where $C$ is a positive constant and
$$\begin{array}{ll}
\displaystyle
J_0(\tau) &
\displaystyle
= \int_{\partial D}\frac{\partial v}{\partial\nu}v\,dS\\
\\
\displaystyle
&
\displaystyle
=\int_D(\vert\nabla v\vert^2+\tau^2\vert v\vert^2)\,dx.
\end{array}
$$

\endproclaim
For the proof of (3.10) see Appendix.  Let us go ahead.
A combination of (3.9) and (3.10) gives
$$\displaystyle
\vert E_1(\tau)^{1/2}-C\tau^{-1/2}J_0(\tau)^{1/2}\vert^2
\le\frac{1}{\tau}\int_{\partial D}\frac{2}{\gamma}
\vert v_1\vert^2\,dx
+C^2\tau^{-1}J_0(\tau)
+O(e^{-2\tau T}).
$$
Since
$$\displaystyle
\vert E_1(\tau)^{1/2}-C\tau^{-1/2}J_0(\tau)^{1/2}\vert^2+C^2\tau^{-1}J_0(\tau)
\ge \frac{1}{2}E_1(\tau),
$$
we finally obtain
$$\begin{array}{ll}
\displaystyle
E_1(\tau)
&
\displaystyle
\le
\frac{1}{\tau}\int_{\partial D}\frac{4}{\gamma}\vert v_1\vert^2\,dx
+4C^2\tau^{-1}J_0(\tau)+O(e^{-2\tau T})\\
\\
\displaystyle
&
\displaystyle
=O(\tau^{-1})
\int_{\partial D}\left(\vert v_1\vert^2+\frac{\partial v}{\partial\nu}v\right)\,dS
+O(e^{-2\tau T})
\end{array}
$$
and hence
$$\displaystyle
\tau^2\int_{\partial D}\vert R_1\vert^2dS
\le C\int_{\partial D}\left(\vert v_1\vert^2+\frac{\partial v}{\partial\nu}v\right)\,dS
+O(\tau e^{-2\tau T}).
\tag {3.11}
$$
Note that we have made use of the fact that $\gamma$ has a positive lower bound.

From (2.8) and (2.11) we have
$$\begin{array}{l}
\displaystyle
\,\,\,\,\,\,
\lim_{\tau\longrightarrow\infty}
e^{2\tau d_{\partial D}(p)}\int_{\partial D}\frac{\partial\tilde{v}}{\partial\nu}\tilde{v}\,dS
\\
\\
\displaystyle
=\frac{\pi}{d_{\partial D}(p)^2}\sum_{q\in\Lambda_{\partial D}(p)}
\frac{1}
{\displaystyle
\sqrt{\text{det}\,(S_q(\partial B_{d_{\partial D}(p)}(p))-S_q(\partial D))}}.
\end{array}
\tag {3.12}
$$
Since $(p-x)/\vert x-p\vert=1$ on $\Lambda_{\partial D}(p)$, it follows from (2.8) and (2.12) for $\gamma\equiv 1$ that
$$\displaystyle
\lim_{\tau\longrightarrow\infty}\tau^{-1}e^{2\tau d_{\partial D}(p)}
\int_{\partial\Omega}\left\vert\frac{\partial\tilde{v}}{\partial\nu}
-\tau\tilde{v}
\right\vert^2\,dS=0.
\tag {3.13}
$$
Since we have (2.9) and (3.3), it follows from (3.11) that
$$\begin{array}{l}
\displaystyle
\,\,\,\,\,\,
\left(\frac{\tau^3}{\varphi(\tau\eta)}\right)^2\tau
\int_{\partial D}\vert R_1\vert^2 dS
\\
\\
\displaystyle
\le
C\left(\sup_{x\in\partial D}\frac{2\gamma}{1+\gamma}\tau^{-1}\int_{\partial D}
\left\vert
\frac{\partial\tilde{v}}{\partial\nu}-\tau\tilde{v}
\right\vert^2dS
+\tau^{-1}\int_{\partial D}\frac{\partial\tilde{v}}{\partial\nu}\tilde{v}dS\right)
+\left(\frac{\tau^3}{\varphi(\tau\eta)}\right)^2O(\tau e^{-2\tau T}).
\end{array}
$$
Now applying (3.12) and (3.13) to this right-hand side and then, using (2.10),
we obtain
$$\displaystyle
\lim_{\tau\longrightarrow\infty}\tau^5e^{2\tau\,\text{dist}\,(D,B)}
\Vert R_1\Vert_{L^2(\partial D)}^2=0.
\tag {3.14}
$$

Now we study the asymptotic behaviour of the second term on (3.6).

\noindent
We have
$$\begin{array}{l}
\displaystyle
\,\,\,\,\,\,
\left\vert
\frac{\displaystyle
\int_{\partial D}\left(\frac{\partial v}{\partial\nu}-\tau\gamma v\right)\,R_1 dS}
{\displaystyle\int_{\partial D}\left(\frac{\partial v}{\partial\nu}-\tau\gamma v\right)\frac{1-\gamma}{1+\gamma}v \,dS}
\right\vert
\\
\\
\displaystyle
\le
\frac{\displaystyle
\tau^{3/2}e^{\tau\text{dist}\,(D,B)}
\left\Vert
\frac{\partial v}{\partial\nu}-\tau\gamma v\right\Vert_{\partial D}
\cdot
\tau^{5/2}e^{\tau\,\text{dist}\,(D,B)}
\Vert R_1\Vert_{L^2(\partial D)}}
{\displaystyle
\tau^4
e^{2\tau\,\text{dist}\,(D,B)}
\left\vert\int_{\partial D}\left(\frac{\partial v}{\partial\nu}-\tau\gamma v\right)\frac{1-\gamma}{1+\gamma}v \,dS
\right\vert
}.
\\
\\
\displaystyle
\end{array}
\tag {3.15}
$$
Again from (2.8) and (2.12) we obtain
$$\begin{array}{l}
\displaystyle
\,\,\,\,\,\,
\lim_{\tau\longrightarrow\infty}
\tau^{-1}e^{2\tau d_{\partial D}(p)}
\int_{\partial D}
\left\vert\frac{\partial\tilde{v}}{\partial\nu}-\tau\gamma\tilde{v}\right\vert^2\,dS
\\
\\
\displaystyle
=\frac{\pi}{d_{\partial D}(p)^2}
\sum_{q\in\Lambda_{\partial D}(p)}
\frac{(1-\gamma(q))^2}
{\displaystyle
\sqrt{\text{det}\,S_q(\partial B_{d_{\partial D}(p)}(p)-S_q(\partial D)}}
\end{array}
$$
and hence (2.10) gives
$$\begin{array}{l}
\displaystyle
\,\,\,\,\,\,
\lim_{\tau\longrightarrow\infty}\tau^{3}e^{2\tau\,\text{dist}\,(D,B)}
\left\Vert
\frac{\partial v}{\partial\nu}-\tau\gamma v
\right\Vert_{L^2(\partial D)}^2
\\
\\
\displaystyle
=\frac{\pi}{4}
\left(\frac{\eta}{d_{\partial D}(p)}\right)^2
\sum_{q\in\Lambda_{\partial D}(p)}
\frac{(1-\gamma(q))^2}
{\displaystyle
\sqrt{\text{det}\,(S_q(\partial B_{d_{\partial D}(p)}(p)-S_q(\partial D))}}.
\end{array}
$$
Applying this, (2.15) and (3.14) to (3.15), we obtain
$$
\lim_{\tau\longrightarrow\infty}
\frac{\displaystyle
\int_{\partial D}\left(\frac{\partial v}{\partial\nu}-\tau\gamma v\right)\,R_1 dS
}
{\displaystyle\int_{\partial D}\left(\frac{\partial v}{\partial\nu}-\tau\gamma v\right)\frac{1-\gamma}{1+\gamma}v \,dS}
=0.
\tag {3.16}
$$
Note also that, if $T$ satisfies (2.6), then, as $\tau\longrightarrow\infty$
$$\begin{array}{l}
\displaystyle
\,\,\,\,\,\,
\tau^4 e^{2\tau\text{dist}\,(D,B)}
\cdot
\tau e^{-\tau T}
(e^{-\tau\,\text{dist}\,(D,B)}+e^{-\tau T})\\
\\
\displaystyle
=\tau^{5}e^{-\tau\,(T-\text{dist}\,(D,B))}
+\tau^5 e^{-2\tau\,(T-\text{dist}\,(D,B))}
\longrightarrow 0.
\end{array}
$$
Now applying this, (2.15), (3.16) to (3.6), we obtain (2.7).

\section{Conclusions and further problems}

In this paper we gave a remark on the previous application
\cite{IEO2} of the enclosure method to an inverse obstacle problem
using the wave governed by the wave equation. We believe that the
argument developed in this paper is quite simple and shall be a
{\it prototype} for applications to other inverse obstacle
problems.

By the way, recently the enclosure method in the time domain has been applied also to the Maxwell system \cite{IMax, IMax2}.
In particular, in \cite{IMax2} we have already obtained a result corresponding to Theorem 1.1.
It is assumed that the electromagnetic field as the solution of the Maxwell system
satisfies the Leontovich boundary condition on the surface of an unknown obstacle.
The Leontovich boundary condition is described by a single positive function defined on the surface.
Thus, it would be interesting to find a formula for the function similar to (1.11) in Theorem 1.2.
This belongs to our next project.

And also developing several applications of the enclosure method in the time domain
to transmission problems for other wave equations and systems should be expected.
See \cite{ISUR} for a survey on recent results obtained by
using the enclosure method.

$$\quad$$

\centerline{{\bf Acknowledgment}}

The author was partially supported by Grant-in-Aid for
Scientific Research (C)(No. 25400155) of Japan  Society for
the Promotion of Science.

$$\quad$$

\section{Appendix. Proof of Lemma 3.1}

This is an application of a reflection argument developed in \cite{LP}.
First of all we compute $(\triangle-\tau^2)R_0$.
Define
$$\displaystyle
\tilde{\phi_{\delta}}(x)
=\frac{1-\tilde{\gamma}(x)}{1+\tilde{\gamma}(x)}\,\phi_{\delta}(x).
$$
Using (4.15) in \cite{LP} (see also \cite{IEE})), we have
$$\begin{array}{ll}
\displaystyle
(\triangle-\tau^2)R_0
&
\displaystyle
=\tilde{\phi}_{\delta}(x)\sum_{i,j}d_{\partial D}(x)a_{ij}(x)(\partial_i\partial_j v)(x^r)
\\
\\
\displaystyle
&
\displaystyle
\,\,\,
+\sum_j\left(\sum_kb_{jk}(x)(\partial_{x_k}\tilde{\phi_{\delta}})(x)
+d_j(x)\tilde{\phi_{\delta}}(x)\right)
(\partial_jv)(x^r)\\
\\
\displaystyle
&
\displaystyle
\,\,\,+
(\triangle\tilde{\phi_{\delta}})(x)v(x^r),
\end{array}
$$
where $a_{ij}(x)$, $b_{jk}(x)$ and $d_j(x)$ with $i,j,k=1,2,3$ are independent of $\tilde{\phi_{\delta}}(x)$,
$v$ and $\tau$; $a_{ij}(x)$ and $b_{jk}(x)$ are $C^1$ and $d_j(x)$ is $C^0$ for $x\in\Bbb R^3\setminus D$
with $d_{\partial D}(x)<2\delta_0$.

Thus we have to study three integrals:
$$\displaystyle
I=\sum_{i,j}\int_{\Bbb R^3\setminus\overline D}
\tilde{\phi_{\delta}}(x)d_{\partial D}(x)a_{ij}(x)(\partial_i\partial_j v)(x^r)\cdot R_1(x)dx,
$$
$$\begin{array}{ll}
\displaystyle
J
&
\displaystyle
=\sum_{j,k}\int_{\Bbb R^3\setminus\overline D}b_{jk}(x)(\partial_k\tilde{\phi_{\delta}})(x)
(\partial_jv)(x^r)\cdot R_1(x)dx\\
\\
\displaystyle
&
\displaystyle
\,\,\,
+\sum_j\int_{\Bbb R^3\setminus\overline D}d_j(x)\tilde{\phi_{\delta}}(x)(\partial_jv)(x^r)R_1(x)dx.
\end{array}
$$
and
$$\displaystyle
K=\int_{\Bbb R^3\setminus\overline D}(\triangle\tilde{\phi_{\delta}})(x)v(x^r)\cdot R_1(x)dx.
$$

In what follows, the symbol $C$ denotes several constants
independent of $\delta$ and $\tau$.

\subsection{Estimating $I$}

The change of variables $x^r=y$ yields
$$\begin{array}{c}
\displaystyle
I=\sum_{i,j}\int_{D}
\tilde{\phi_{\delta}}(y^r)d_{\partial D}(y^r)a_{ij}(y^r)(\partial_i\partial_j v)(y)\cdot R_1(y^r)J(y)dy,
\end{array}
$$
where $J(y)$ denotes the Jacobian.
Since $d_{\partial D}(y^r)=0$ on $\partial D$, integration by parts yields
$$
\displaystyle
I
=-\sum_{i,j}\int_{D}
\partial_i(\tilde{\phi_{\delta}}(y^r)d_{\partial D}(y^r)a_{ij}(y^r)R_1(y^r)J(y))(\partial_j v)(y)dy.
$$
Set
$$\displaystyle
R_1^r(y)=R_1(y^r)
$$
and
$$\displaystyle
D_{\delta}=\{y\in D\,\vert\,d_{\partial D}(y)<2\delta\}.
$$
Noting $d_{\partial D}(y^r)=d_{\partial D}(y)$, we can easily obtain
$$\displaystyle
\vert I\vert
\le C(\Vert R_1^r\Vert_{L^2(D_{\delta})}
+\delta\Vert\nabla R_1^r\Vert_{L^2(D_{\delta})})\Vert\nabla v\Vert_{L^2(D_{\delta})}.
\tag {A.1}
$$

\subsection{Estimating $J$}

Since
$$\begin{array}{ll}
\displaystyle
J
&
\displaystyle
=\sum_{j,k}\int_{D}b_{jk}(y^r)(\partial_k\tilde{\phi_{\delta}})(y^r)\partial_jv(y)\cdot R_1(y^r)J(y)dy\\
\\
\displaystyle
&
\displaystyle
\,\,\,
+\sum_j\int_{D}d_j(y^r)\tilde{\phi_{\delta}}(y^r)\partial_jv(y)R_1(y^r)J(y)dy,
\end{array}
$$
we have
$$\displaystyle
\vert J\vert
\le C\delta^{-1}\Vert R_1^r\Vert_{L^2(D_{\delta})}\Vert\nabla v\Vert_{L^2(D_{\delta})}.
\tag {A.2}
$$

\subsection{Estimating $K$}

Since
$$\displaystyle
K=\int_{D}(\triangle\tilde{\phi_{\delta}})(y^r)v(y)\cdot R_1(y^r)J(y)dy,
$$
we have
$$\displaystyle
\vert K\vert
\le C\delta^{-2}\Vert R_1^r\Vert_{L^2(D_{\delta})}\Vert v\Vert_{L^2(D_{\delta})}.
\tag {A.3}
$$

Summing up (A.1)-(A.3), we obtain
$$\begin{array}{ll}
\displaystyle
\left\vert\int_{\Bbb R^3\setminus\overline D}(\triangle-\tau^2)R_0\cdot R_1 dx\right\vert
&
\displaystyle
\le
C\left(\delta\Vert\nabla R_1^r\Vert_{L^2(D_{\delta})}\Vert\nabla v\Vert_{L^2(D)}\right.
\\
\\
\displaystyle
&
\displaystyle
\,\,\,
+\delta^{-1}\Vert R_1^r\Vert_{L^2(D_{\delta})}\Vert\nabla v\Vert_{L^2(D)}
\\
\\
\displaystyle
&
\displaystyle
\,\,\,
\left.
+\delta^{-2}\Vert R_1^r\Vert_{L^2(D_{\delta})}\Vert v\Vert_{L^2(D)}\right).
\end{array}
\tag {A.4}
$$

Here note that
$$
\left\{
\begin{array}{l}
\displaystyle
\Vert\nabla v\Vert_{L^2(D)}\le J_0(\tau)^{1/2},
\\
\\
\displaystyle
\Vert v\Vert_{L^2(D)}\le\tau^{-1}J_0(\tau)^{1/2}
\end{array}
\right.
\tag {A.5}
$$
and
$$\left\{
\begin{array}{l}
\displaystyle
\Vert\nabla R_1^r\Vert_{L^2(D_{\delta})}
\le C
E_1^0(\tau)^{1/2},
\\
\\
\displaystyle
\Vert R_1^r\Vert_{L^2(D_{\delta})}\le\tau^{-1}E_1^0(\tau)^{1/2},
\end{array}
\right.
\tag {A.6}
$$
where
$$\displaystyle
E_1^0(\tau)
=\int_{\Bbb R^3\setminus\overline D}(\vert\nabla R\vert^2+\tau^2\vert R\vert^2)\,dx.
$$

Applying (A.5) and (A.6) to the right-hand side of (A.4), we obtain
$$\begin{array}{l}
\displaystyle
\,\,\,\,\,\,
\left\vert\int_{\Bbb R^3\setminus\overline D}(\triangle-\tau^2)\tilde{v}\cdot\epsilon dx\right\vert
\\
\\
\displaystyle
\le
C(\delta
+\delta^{-1}\tau^{-1}
+\delta^{-2}\tau^{-2})(J_0(\tau)E_1^0(\tau))^{1/2}.
\end{array}
$$
Therefore, choosing $\delta=\tau^{-1/2}$, we obtain
$$\begin{array}{c}
\displaystyle
\left\vert\int_{\Bbb R^3\setminus\overline D}(\triangle-\tau^2)R_0\cdot R_1 dx\right\vert
\le
C\tau^{-1/2}(J_0(\tau)E_1^0(\tau))^{1/2}.
\end{array}
$$
Since $E_1(\tau)\ge E_1^0(\tau)$, this yields the desired estimate.

\vskip1cm
\noindent
e-mail address

ikehata@amath.hiroshima-u.ac.jp


\begin{thebibliography}{99}




\bibitem{BH}  Bleistein, N. and Handelsman, R. A.,
                Asymptotic expansions of integrals, New York, Dover, 1986.






\bibitem{CH}  Courant, R. and Hilbert, D.,
              Methoden der Mathematischen Physik, vol. 2., Berlin, Springer, 1937.







\bibitem{DL}  Dautray, R. and Lions, J-L., Mathematical analysis and numerical methods for
              sciences and technology,
              Evolution problems I, Vol. {\bf 5}, Springer-Verlag, Berlin, 1992.










\bibitem{GT}  Gilbarg, D. and Trudinger, N. S., Elliptic partial differential equations
              of second order, second.ed. (Berlin:Springer), 1983.










\bibitem{I1} Ikehata, M.,
             \newblock Enclosing a polygonal cavity in a two-dimensional bounded domain from Cauchy data,
             \newblock Inverse Problems, {\bf 15}(1999), 1231-1241.







\bibitem{IEO2} Ikehata, M.,
              The enclosure method for inverse obstacle scattering problems with dynamical data over a finite time
              interval: II. Obstacles with a dissipative boundary or finite refractive index and back-scattering data,
              Inverse Problems, {\bf 28}(2012) 045010 (29pp).






\bibitem{IEE} Ikehata, M.,
              Extracting the geometry of an obstacle and a zeroth-order coefficient
              of a boundary condition via the enclosure method using a single reflected wave over a finite
              time interval,
              Inverse Problems, {\bf 30}(2014) 045011 (24pp).







\bibitem{ISUR} Ikehata, M.,
             New development of the enclosure method for inverse obstacle scattering,
             to appear as Capter  6, {\it Inverse Problems and Computational Mechanics}
             (eds. L. Marin, L. Munteanu, V. Chiroiu),
             Vol. {\bf 2} , Editura Academiei, Bucharest, Romania, 2015.






\bibitem{IMax}  Ikehata, M.,
             The enclosure method for inverse obstacle scattering using a single electromagnetic
             wave in time domain,  Inverse Problems and Imaging, to appear.





\bibitem{IMax2} Ikehata, M.,
             On finding an obstacle with the Leontovich boundary condition via the time domain
             enclosure method, arXiv:1510.08209v1 [math.AP] 28 Oct 2015.









\bibitem{LP}  Lax, P. D. and Phillips, R. S., The scattering of sound waves by an obstacle,
              Comm. Pure and Appl. Math., {\bf 30}(1977), 195-233.









\bibitem{Mo} Majda, A.,
             High frequency asymptotics for the scattering matrix and the inverse problem of acoustic scattering,
             Comm. Pure and Appl. Math., {\bf 29}(1976), 261-291.








\bibitem{O}  O'Neill, B.,  Elementary Differential Geometry, Revised , 2nd Edition, 2006, Academic Press.













\end{thebibliography}
\end{document}